\documentclass{amsart}
\usepackage[all]{xy}
\usepackage{amsfonts,amscd,amssymb,amsmath,mathrsfs}
\newtheorem{thm}{Theorem}[section]
\newtheorem{lem}[thm]{Lemma}
\newtheorem{prop}[thm]{Proposition}
\newtheorem{cor}[thm]{Corollary}
\theoremstyle{definition}
\newtheorem{defn}{Definition}[section]

\theoremstyle{remark}

\numberwithin{equation}{section}

\begin{document}

\newcommand{\Soc}{\operatorname{Soc}}
\newcommand{\modP}{\mod_{P}\Lambda}
\newcommand{\modl}{\mod \Lambda}
\newcommand{\module}{\operatorname{mod}}
\newcommand{\p}{\operatorname{p}}
\newcommand{\inj}{\operatorname{i}}
\newcommand{\Module}{\operatorname{Mod}}
\newcommand{\Cok}{\operatorname{Coker}}
\newcommand{\Hom}{\operatorname{Hom}}
\newcommand{\h}{\operatorname{h}}
\newcommand{\e}{\operatorname{e}}
\newcommand{\res}{\operatorname{res}}
\newcommand{\Tr}{\operatorname{Tr}}
\newcommand{\TrD}{\operatorname{TrD}}
\newcommand{\rad}{\operatorname {{\bold r}}}
\newcommand{\La}{\operatorname{\Lambda}}
\newcommand\End{\operatorname{End}}
\newcommand\Ext{\operatorname{Ext^1_\Lambda}}
\newcommand\Ex{\operatorname{Ext}}
\newcommand\ann{\operatorname{ann}}
\newcommand\coend{\operatorname{coend}}
\newcommand\Img{\operatorname{Im}}
\newcommand\D{\operatorname{D}}
\newcommand\DTr{\operatorname{DTr}}
\newcommand\Ker{\operatorname{Ker}}
\newcommand\Coker{\operatorname{Coker}}
\newcommand{\tr}{\operatorname {t}}
\newcommand{\ct}{\operatorname{ct}}
\newcommand{\rej}{\operatorname{rej}}
\newcommand{\gd}{\operatorname{{gl.dim}}}
\newcommand{\radic}{\operatorname{rad}}
\newcommand{\add}{\operatorname{add}}
\newcommand{\Ind}{\operatorname{Ind}}
\newcommand{\Sub}{\operatorname{Sub}}
\newcommand{\Fac}{\operatorname{Fac}}
\newcommand{\lap}{\operatorname{l}}
\newcommand{\rap}{\operatorname{r}}
\newcommand{\locn}{\operatorname{lnRep}}
\newcommand{\rep}{\operatorname{Rep}}

\newcommand{\spec}{(\Gamma,\Lambda,\textrm{\textbf{d}},\mathfrak B)}
\newcommand{\newspec}{(\Gamma,\sigma_x\Lambda)}
\newcommand{\gl}{(\Gamma,\Lambda)}
\newcommand{\poset}{(\Gamma_{0},\Lambda)}
\newcommand{\Stilde}{\tilde S}
\newcommand{\Ttilde}{\tilde T}
\newcommand{\ad}{(+){\rm -admissible}}
\newcommand{\pa}{{\tilde{\mathcal P}}(A)}
\newcommand{\modA}{{\rm mod}A}
\newcommand{\hull}{H_{\Lambda}}
\newcommand{\modr}{\mathrm{mod}\,R_{\gl}}
\newcommand{\ds}[1]{\ensuremath{#1}}
\newcommand{\quiver}[2]{\ensuremath{\left( #1, #2 \right) }}
\newcommand{\quiverA}[1]{\ensuremath{\left( A_n, #1 \right) }}
\newcommand{\set}[1]{\ensuremath{\left\{#1\right\}}}

\bibliographystyle{plain}
\title[Valued quivers and the Weyl group]{Preprojective representations of valued quivers
and reduced words in the Weyl group of a Kac-Moody algebra}


\author{Mark Kleiner and Allen Pelley}
\address{Department of Mathematics, Syracuse University, Syracuse, New
York 13244-1150}
\email{mkleiner@syr.edu}
\email{anpelley@syr.edu}
\thanks{The authors are supported by the NSA grant H98230-06-1-0043.}

\date{}

\dedicatory{Dedicated to the memory of Andrei Vladimirovich Roiter}

\commby{}

\begin{abstract}This paper studies connections between the preprojective representations of a valued quiver, the (+)-admissible sequences of vertices, and the Weyl group by associating to each preprojective   representation a canonical (+)-admissible sequence.  A (+)-admissible sequence is the canonical sequence of some preprojective representation if and only if the product of  simple reflections associated to the vertices of the sequence  is a reduced word in the Weyl group.  As a consequence, for any Coxeter element of the Weyl group associated to an indecomposable symmetrizable generalized Cartan matrix, the group is infinite if and only if the powers of the element are reduced words.  The latter strengthens known results of Howlett, Fomin-Zelevinsky, and the authors. 
\end{abstract}

\maketitle

\section*{Introduction}\label{Section:intro}

Let $\mathcal W$ be a Coxeter group generated by  reflections $\sigma_1,\dots,\sigma_n$ ~\cite{bourbLie4-6}, and let  $c$ be any Coxeter element of  $\mathcal W$, i.e., $c=\sigma_{x_n}\dots\sigma_{x_1}$  where $x_1,\dots,x_n$ is any permutation of the  numbers $1,\dots,n$.  Recall that an indecomposable symmetrizable generalized
Cartan matrix (see ~\cite{kac1990}) is an
integral  $n\times n$  matrix $A=(a_{ij})$ with  $a_{ii} = 2$ and $a_{ij}\le0$ for $i\ne j$ that
is not conjugate under a permutation matrix to a block-diagonal
matrix $\left(\begin{matrix}A_1&0\\0&A_2\end{matrix}\right)$ and such that
there exist  integers $d_i>0$ satisfying $d_i a_{ij}=d_j a_{ji}$ for all $i,j$.
We identify the root lattice $Q$ associated with $A$ with the free
abelian group ${\mathbb Z}^n$ by identifying the simple roots
$\alpha_1,\dots,\alpha_n$ of $Q$ with the standard basis vectors
$e_1,\dots,e_n$ of ${\mathbb Z}^n$.  Then the simple
reflections on $Q$ identify with the reflections
$\sigma_1,\dots,\sigma_n$ on ${\mathbb Z}^n$ given by
$\sigma_i(e_j)=e_j-a_{ij}e_i$ for all $i,j$, and the Weyl group
$\mathcal W(A)$ is the subgroup of $GL({\mathbb Z}^n)$ generated by
$\sigma_1,\dots,\sigma_n$.  Andrei Zelevinsky brought to our attention the following two results.   Howlett proved that  $\mathcal W$ is  infinite if and only if $c$ has infinite order ~\cite[Theorem 4.1]{h}.  Fomin and Zelevinsky proved the  following. Let $A$ be bipartite, i.e., the set $\{1,\dots,n\}$ is a disjoint union of nonempty subsets $I,J$ and, for $h\ne l$,  $a_{hl}=0$ if either $h,l\in I$ or $h,l\in J$.  If $c=\prod_{i\in I}\sigma_i\prod_{j\in J}\sigma_j$, then ${\mathcal W}(A)$ is infinite if and only if the powers of $c$ are reduced words in the $\sigma_h$'s ~\cite[Corollary 9.6]{CAIV}.  Inspired by the latter, we proved ~\cite[Theorem 4.8]{kp} that if  $A$ is symmetric and $c$ is any Coxeter element, then ${\mathcal W}(A)$ is infinite if and only if the powers of $c$ are reduced words, which strengthens the aforementioned results of Howlett and Fomin-Zelevinsky in the case when ${\mathcal W}={\mathcal W}(A)$ and $A$ is symmetric.  In this paper we obtain a further strengthening (Theorem \ref{infWeylGrp}) by removing the additional assumption that $A$ is symmetric.

The proof of the indicated result from ~\cite{kp} was based on the interplay  between the category  $\tilde{\mathscr P}$ of preprojective modules ((+)-irregular representations) over the path algebra of a finite connected quiver without oriented cycles,  the set $\mathfrak S$ of (+)-admissible sequences, and  the Weyl group $\mathcal W$ of the underlying (nonoriented) graph of the quiver, as defined by Bernstein, Gelfand, and Ponomarev ~\cite{bgp}.  We used the  one-to-one correspondence between the finite connected graphs (with no loops but multiple edges allowed) and the indecomposable symmetric Cartan matrices, and relied on the study of $\tilde{\mathscr P}$ in terms of the combinatorics of $\mathfrak S$ undertaken in ~\cite{kt}.  Since there is a one-to-one correspondence between the connected valued graphs and the indecomposable symmetrizable Cartan matrices ~\cite[p. 1]{dr}, we replace graphs with valued graphs, representations of quivers with representations of valued quivers studied in ~\cite{dr}, and rely on the extension of the results of ~\cite{kt} to representations of valued quivers provided in ~\cite{kt2}.  Our study of $\tilde{\mathscr P}$,  $\mathfrak S$, and $\mathcal W$ for a valued quiver is also of independent interest.   We associate  to any $M\in\tilde{\mathscr P}$ a canonical sequence $S_M\in\mathfrak S$ and consider, for each $S\in\mathfrak S$,  the element $w(S)\in \mathcal W$ that is the composition of  simple reflections associated to the vertices of $S$.  Then if $M,N\in\tilde{\mathscr P}$ are indecomposable, we  have  $M\cong N$ if and only if $w(S_M)=w(S_N)$ (Theorem \ref{reducedword}(b)).  For all $S\in\mathfrak S$, there exists an $M\in\tilde{\mathscr P}$ satisfying $S\sim S_M$ if and only if the word $w(S)\in \mathcal W$ is reduced (Theorems \ref{princseqreducedwrd} and \ref{seqreducedwrd}).  If the valued graph is not a Dynkin diagram of the type $A_n,B_n,C_n, D_n,  E_6,E_7,E_8,F_4,$ or $G_2,$ then for all $S\in\mathfrak S$, the word $w(S)$ is reduced and there exists an $M\in\tilde{\mathscr P}$ satisfying $S\sim S_M$ (Corollary \ref{reduced}).  We note that the connection between   Coxeter-sortable elements of $\mathcal W$ and  the canonical sequences in $\mathfrak S$ associated to  the modules in $\tilde{\mathscr P}$, as explained in  ~\cite[Section 5]{kp},  holds for valued quivers.  Nathan Reading introduced Coxeter-sortable elements for an arbitrary Coxeter group ~\cite{rea} and proved that if the group is finite, the set of Coxeter-sortable elements maps bijectively onto  the set of clusters ~\cite{FZ2003, mrz2003} and onto the set of noncrossing partitions ~\cite{McCam}.  It follows that to each preprojective representation of a valued quiver correspond a cluster and a noncrossing partition, and we will study these clusters and noncrossing partitions in the future.

We leave it to the reader to dualize the results of this paper by replacing preprojective representations of valued quivers with preinjective ones, and (+)-admissible sequences with ($-$)-admissible sequences ~\cite{dr}.  

\section{Preliminaries}\label{prelim}

We recall some facts, definitions, and notation, using freely
~\cite{ars,bgp,dr, kt2}. A {\it graph} is a pair
$\Gamma=(\Gamma_{0},\Gamma_{1})$, where $\Gamma_{0}$ is the set of
vertices, and  the set of edges $\Gamma_{1}$ consists of some two-element
subsets of $\Gamma_0$.  A {\it valuation}
$\textrm{\textbf{b}}$ of  $\Gamma$ is a set of 
integers $b_{ij}\ge0$ for all pairs $i,j\in\Gamma_0$ where
$b_{ii}=0$; if $i\ne j$ then $b_{ij}\ne0$ if and only if $\{i,j\}\in\Gamma_1$; and there exist integers $d_i>0$ satisfying
$d_ib_{ij}=d_jb_{ji}\ \textrm{for\ all}\ i, j\in\Gamma_0.$
The pair $(\Gamma,\textrm{\textbf{b}})$ is a {\it valued graph}, which is called
{\it connected\,} if   $\Gamma$ is connected.  

An orientation, $\Lambda$, on $\Gamma$ consists of two functions
$s:\Gamma_{1}\to\Gamma_{0}$ and $e:\Gamma_{1}\to\Gamma_{0}$.  If $a\in\Gamma_{1}$, then $s(a)$ and $e(a)$ are the vertices incident
with $a$, called the starting point and the endpoint of
$a$, respectively.  The  triple
$(\Gamma,\textrm{\textbf{b}},\Lambda)$ is  a {\it valued
quiver}, and $a$ is then called an arrow of the quiver.  Given a sequence of arrows
$a_{1},\dots,a_{t},\ t>0,$ satisfying $e(a_{i})=s(a_{i+1}),\
0<i<t,$ one forms a path $p=a_{t}¥\dots a_{1}$ of length $t$ in
$(\Gamma,\textrm{\textbf{b}},\Lambda)$ with
$s(p)=s(a_{1})$ and $e(p)=e(a_{t})$. By definition,
for all $x\in\Gamma_{0}$ there is a unique path $e_x$ of length $0$ with
$s(e_x)=e(e_x)=x$. We say that $p$ is a path from $s(p)$ to $e(p)$  and write $p:s(p)\to e(p)$.  A path $p$ of length $>0$ is an
oriented cycle if $s(p)=e(p)$. The set of vertices of any valued
quiver without oriented cycles (no finiteness assumptions) becomes a partially ordered set (poset) if one puts $x\le y$
whenever there is a path from $x$ to $y$.  If
$(\Gamma,\textrm{\textbf{b}},\Lambda)$ has no oriented cycles, we
denote this poset by $(\Gamma_{0}¥,\La)$.  Throughout, all orientations
$\Lambda,\Theta$, etc., are such that
$(\Gamma,\textrm{\textbf{b}},\Lambda)$,
$(\Gamma,\textrm{\textbf{b}},\Theta)$, etc., have no oriented
cycles.

To define representations of a valued quiver
$(\Gamma,\textrm{\textbf{b}},\Lambda)$, one has to choose a {\it
modulation} $\mathfrak B$ of the valued graph
$(\Gamma,\textrm{\textbf{b}})$, which by definition is a set of
division rings $\mathbf{k}_i$, $i\in\Gamma_0$, together with a
$\mathbf{k}_i-\mathbf{k}_j$-bimodule $_iB_j$ and a
$\mathbf{k}_j-\mathbf{k}_i$-bimodule $_jB_i$ for each
$\{i,j\}\in\Gamma_1$ such that

(i) there are $\mathbf{k}_j-\mathbf{k}_i$-bimodule isomorphisms\vskip.05in
\centerline{$_jB_i\cong\textrm{Hom}_{\mathbf{k}_i}(_iB_j,\mathbf{k}_i)\cong\textrm{Hom}_{\mathbf{k}_j}(_iB_j,\mathbf{k}_j)$} \noindent and

(ii) dim$_{\mathbf{k}_i}({_iB_j})=b_{ij}$.

Unless indicated otherwise, for the rest of the paper we  denote by  $\Gamma$ an arbitrary finite connected valued graph with $|\Gamma_{0}|>1$, where $|X|$ stands for the cardinality
of a set $X$, and with
a valuation $\textrm{\textbf{b}}$ and modulation $\mathfrak
B$; denote by $\gl$ the  corresponding valued quiver with
orientation $\La$; and assume that $\dim_k\mathbf{k}_i<\infty$ for all $i$,
where $k$ is a common central subfield  of the $\mathbf{k}_i$'s
acting centrally on all bimodules $_iB_j$.  Under the assumption,
each   $_iB_j$  is a finite dimensional $k$-space.

A (left) representation $(V,f)$ of $\gl$ is a set of finite
dimensional left $\mathbf{k}_i$-spaces $V_i$, $i\in\Gamma_0$,
together with $\mathbf{k}_j$-linear maps
$f_a:{_jB_i}\otimes_{\mathbf{k}_i}V_i\to V_j$ for all arrows
$a:i\to j$.   Defining morphisms of representations in a
natural way, one obtains the category $\textrm{Rep}\gl$ of
representations of the valued quiver $\gl$.

Putting $\mathbf{k}=\prod_{i\in\Gamma_0}\mathbf{k}_i$ and viewing
$B=\underset{i\to j}\bigoplus{_jB_i}$ as a
$\mathbf{k}$-$\mathbf{k}$-bimodule where $\mathbf{k}$ acts on
${_jB_i}$ from the left via the projection
$\mathbf{k}\to\mathbf{k}_j$ and from the right via the projection
$\mathbf{k}\to\mathbf{k}_i$, one forms the tensor ring
$\textrm{T}(\mathbf{k},B)=\bigoplus_{n=0}^\infty B^{(n)}$ where
$B^{(n)}=B\otimes_\mathbf{k}\cdots\otimes_\mathbf{k} B$ is the
$n$-fold tensor product, and the multiplication is given by the
isomorphisms $B^{(n)}\otimes B^{(m)}\to B^{(n+m)}$ ~\cite[p.
386]{dr1}.  Since $\gl$ has no oriented cycles,
$\textrm{T}(\mathbf{k},B)$ is a finite dimensional $k$-algebra and
we denote it by $k\gl$.  Let  $e_i\in\mathbf{k}$ be the $n$-tuple
that has $1\in\mathbf{k}_i$ in the $i$th place and 0 elsewhere.  A
left $k\gl$-module $M$ is {\it finite dimensional} if
$\mathrm{dim}_{\mathbf{k}_i}e_iM<\infty$ for all $i$, which is
equivalent to $\dim_k M<\infty$. We let f.d.$\,k\gl$ denote the
category of finite dimensional left $k\gl$-modules. The categories
Rep$\gl$ and f.d.$\,k\gl$ are equivalent ~\cite[Proposition
10.1]{dr1} and we view the equivalence as an identification.  In
this paper all $k(\Gamma,\Theta)$-modules are finite dimensional, for  any orientation  $\Theta$ without oriented cycles.

If $\gl$ is a valued quiver  and $x\in\Gamma_{0}$, let
$\sigma_{x}\Lambda$ be the orientation on $\Gamma$ obtained by
reversing the direction of each arrow incident with $x$ and
preserving the directions of the remaining arrows.  There results a
new valued quiver $\newspec$ (with the same  valuation $\mathbf b$
and modulation $\mathfrak B$).  A vertex $x$ is a {\it sink}  if no arrow starts at $x$.
For a sink $x$, we recall the definition of the {\it reflection} functor
$F^{+}_{x}:\textrm{Rep}\gl\to\textrm{Rep}\newspec$
~\cite[pp. 15-16]{dr}.

Let $(V,f)\in\textrm{Rep}\gl$ and let $(W,g)=F^{+}_{x}(V,f)$.  Then
$W_y=V_y$ for all $y\neq x$, and $g_b=f_b$ for  the arrows $b$
of $\newspec$ that do not start at $x$.  Let $a_i:y_i\to x,\
i=1,\dots,l$, be the arrows of $\gl$ ending at $x$. Then the
reversed arrows $a_i':x\to y_i$ are the arrows of
$\newspec$ starting at $x$.  Consider the exact sequence
\newline\centerline{$
0\to\Ker
h\overset{j}\to\overset{l}{\underset{i=1}\oplus}{_xB_{y_i}}\otimes_{\mathbf{k}_
{y_i}} V_{y_i}\overset{h}\to V_x
$}
of $\mathbf{k}_x$-spaces where $h$ is induced by the maps
$f_{a_i}:{_xB_{y_i}}\otimes_{\mathbf{k}_ {y_i}}V_{y_i}\to V_x$. Then
$W_x=\Ker h$ and  each map
$g_{a_i'}:{_{y_i}B_x}\otimes_{\mathbf{k}_x}W_x\to W_{y_i}=V_{y_i}$
is obtained from the map
$W_x\to{_xB_{y_i}}\otimes_{\mathbf{k}_{y_i}}W_{y_i}$ induced by $j$
using the isomorphisms  below ~\cite[pp.
14-15]{dr}.\vskip.04in
\centerline{$\textrm{Hom}_{\mathbf{k}_x}(W_x,{_xB_{y_i}}\otimes_{\mathbf{k}_{y_i}}W_{y_i})\cong
\textrm{Hom}_{\mathbf{k}_x}(W_x,\textrm{Hom}_{\mathbf{k}_{y_i}}({_{y_i}B_x},\mathbf{k}_{y_i})\otimes_{\mathbf{k}_{y_i}}W_{y_i})\cong$}\vskip.04in
\centerline{$\textrm{Hom}_{\mathbf{k}_x}(W_x,\textrm{Hom}_{\mathbf{k}_{y_i}}({_{y_i}B_x},W_{y_i}))\cong\textrm{Hom}_{\mathbf{k}_{y_i}}({_{y_i}B_x}\otimes_{\mathbf{k}_{x}}W_x,W_{y_i})$}\vskip.04in

A sequence of vertices $S=x_{1},x_{2},\dots,x_{s},\ s\ge0,$ is
called {\it (+)-admissible} on $\gl$ if it either is  empty, or
satisfies the following conditions: $x_{1}$ is a sink in $\gl$, $x_{2}$ is a sink in
$(\Gamma,\sigma_{x_{1}}\Lambda)$, and so on.  We denote by $\mathfrak S$ the set of (+)-admissible sequences on
$\gl$. 

The following is \cite[Definition 1.1]{kt2} and ~\cite[Definition 2.1]{kp}.

\begin{defn}
If $S=x_1,\dots,x_s,\ s\ge0,$ is in $\mathfrak S$, we write
$\La^S=\sigma_{x_s}\dots\sigma_{x_1}\La$ and, in particular,
$\La^{\emptyset}=\La$. The {\it support} of $S$, $\mathrm{Supp}\,
S$,  is the set of distinct vertices among $x_j$, $1\leq j\leq s$.
The {\it length} of $S$ is
$\ell(S)=s$; the \textit{multiplicity} of $v\in\Gamma_0$ in $S$,
$m_S(v)$, is the (nonnegative) number of subscripts $j$ satisfying
$x_j=v$. A sequence $K\in\mathfrak S$ is \textit{complete} if
$m_K(v)=1$ for all $v\in\Gamma_0$. If $S=x_1,\dots,x_s$ and
$T=y_1,\dots,y_t$ are $(+)$-admissible on $\gl$ and
$(\Gamma,\Lambda^{S})$, respectively, the concatenation of $S$ and
$T$ is the sequence $ST=x_1,\dots,x_s,y_1,\dots,y_t$ in $\mathfrak S$.  If $K$ is
complete, then $\Lambda^K=\Lambda$, so 
the concatenation $K^m$ of $m>0$ copies of $K$ is in $\mathfrak S$.
\end{defn}

We quote \cite[Definition 1.2]{kt2}.

\begin{defn}\label{defequiv}
If  $S=x_1,\dots,x_{i},x_{i+1},\dots,x_s,\ 0<i<s$, is in
$\mathfrak S$ and no edge of $\Gamma$ connects
$x_{i}$ with $x_{i+1}$, then
$T=x_1,\dots,x_{i+1}¥,x_{i}¥,\dots,x_s$ is in $\mathfrak S$ and we
set $SrT$. Let $\sim$ be the equivalence relation that is a
reflexive and transitive closure of the symmetric binary relation
$r$.
\end{defn}

If  $S=x_{1},\dots,x_{s}$ is in $\mathfrak S$, we set
$F(S)=F^{+}_{x_{s}}\dots F^{+}_{x_{1}}:\textrm{Rep}\gl\to
\textrm{Rep}(\Gamma, \Lambda^S)$.  It follows from the analog of  ~\cite[Lemma
1.2, proof of part 3)]{bgp} for representations of valued quivers
that $S\sim T$ implies $F(S)=F(T)$.  If $K\in\mathfrak S$ is complete, then $F(K)=\Phi^{+}$ is  the {\it Coxeter} functor
~\cite[p. 19]{dr}.   We say that $S\in\mathfrak S$
\textit{annihilates} $M\in\mathrm{f.d.}\,k\gl$ if
$F(S)(V,f)=0$ where $(V,f)\in\mathrm{Rep}\gl$ is 
identified with $M$. In light of this identification, we often write
$F(S)M$ or $\Phi^+ M$.

A {\it source} is a vertex of a quiver at which no arrow ends.
Replacing sinks with sources, one gets similar definitions of a {\it
reflection} functor  $F_{x}^{-}$, a ($-$){\it -admissible} sequence,
and the {\it Coxeter} functor $\Phi^{-} $  ~\cite{dr}.

The definition of a
preprojective representation ~\cite[p. 22]{dr} is equivalent to the following
(see \cite[Note 2]{bgp}).

\begin{defn}\label{intr2} A module $M\in\mathrm{f.d.}\,k\gl$ is
{\it preprojective} if there exists an $S\in\mathfrak S$ that
annihilates it.
\end{defn}

We quote \cite{kp,kt,kt2} for the use in
Section \ref{Section:Weylgroup} and note that any result for ordinary quivers that deals with $\mathfrak S$ but not with $\tilde{\mathscr P}$ holds for valued quivers and no proof is needed, for (+)-admissible sequences depend on the graph $\Gamma$ and orientation $\Lambda$, but not on
the valuation $\mathbf b$ or the modulation $\mathfrak B$.

The following is \cite[Proposition
1.9]{kt}.

\begin{prop}\label{Prop:canonical}
Let ${S \in \mathfrak{S}}$ be nonempty.

\begin{itemize}
\item[(a)] ${S \sim S_1 \dots S_r}$ where, for all ${i}$,
${S_{i}}\in\mathfrak{S}$ consists of distinct vertices and 
$\,\mathrm{Supp}\,{S_{i}}=$ $\mathrm{Supp}\,{S_{i}}\dots S_r$.
If ${\mathrm{Supp}\,{S_{i}}\neq\Gamma_{0}}$ then
$\mathrm{Supp}\,{S_{i+1}}\subsetneq\mathrm{Supp}\,{S_{i}}$.

\item[(b)] If ${T\sim T_{1}}\dots T_{q}\,$ 
in ${\mathfrak S}$ is nonempty where, for all ${j}$,  ${T_{j}}\in\mathfrak{S}$ consists of
distinct vertices  and $\,{\mathrm{Supp}\,{T_{j}}=\mathrm{Supp}\,{T_{j}\dots
T_{q}}}$, then ${S\sim T}$ if and only if ${r=q}$ and
${S_{i}\sim T_{i}}$ on ${(\Gamma,\Lambda^{S_{1}\dots
S_{i-1}}),i=1,\dots,r}$.

\end{itemize}
\end{prop}

The sequence ${S_1  \dots S_r}$ of
Proposition \ref{Prop:canonical}(a) is the \textit{canonical
form} of ${S \in \mathfrak{S}}$,  and the integer $r$ is the {\it size} of $S$.  

The following  is ~\cite[Definition
2.1]{kt}.

\begin{defn}
If ${S,T \in \mathfrak{S}}$, we say that ${S}$ is a
{\it subsequence} of ${T}$ and write ${S\preccurlyeq T}$ if
${T\sim SU}$ for some (+)-admissible sequence ${U}$.
\end{defn}

We quote \cite[Proposition 1.5]{kt2} and \cite[Proposition 2.3]{kp}.

\begin{prop}\label{subsequence}\begin{itemize}\item[(a)] The relation $\preccurlyeq$ is a preorder.
\end{itemize}
Let $S,T\in\mathfrak S$.
\begin{itemize}
\item[(b)] We have $S\preccurlyeq T$ and $T\preccurlyeq S$ if and only if
$S\sim T$.
\item[(c)] $S\preccurlyeq T$ if and only if  for all $v\in\Gamma_0$, $m_S(v)\le m_T(v)$.
\end{itemize}
\end{prop}

We quote \cite[Definition 2.4 and Proposition 2.6]{kp}.

\begin{prop}\label{def:meetjoin}  The poset of equivalence classes in $\mathfrak S$ is a
lattice where the greatest lower bound $\wedge$ and the least upper bound $\vee$ are as follows.

Let ${S, T \in \mathfrak{S}}$ be nonempty and let ${S_{1} \dots
S_{r}},\ { T_{1} \dots T_{q}}$ be their canonical forms,
respectively, where ${r
\leq q}$.  Then:
\begin{itemize}
\item[(a)] If $\,\mathrm{Supp}\,{S} \cap
\mathrm{Supp}\,{T}=\emptyset$ then ${S \wedge T}=\emptyset$; and if $\mathrm{Supp}\,{S} \cap
\mathrm{Supp}\,{T}\ne\emptyset$, the canonical form of  ${S \wedge T}$ is  ${R_{1} \dots
R_{s}}$ where $s\le q,r$ is the largest integer satisfying
${\mathrm{Supp}\,{R_{i}} = \mathrm{Supp}\,{S_{i}} \cap
\mathrm{Supp}\,{T_{i}}}\ne\emptyset$ for $0< i \leq s$.

\item[(b)] The
canonical form of ${S \vee T}$ is  ${R_{1} \dots R_{q}}$ where
$\mathrm{Supp}\,{R_{i}} = \mathrm{Supp}\,{S_{i}}\cup
\mathrm{Supp}\,{T_{i}}$ for ${0< i \leq r}$, and
${\mathrm{Supp}\,{R_{i}} = \mathrm{Supp}\,{T_{i}}}$ for ${r < i \leq
q}$.
\end{itemize}

For all
$S\in\mathfrak{S}$, we have  ${S \wedge \emptyset=\emptyset}$  and
${S \vee \emptyset=S}$.
\end{prop}

 Recall  that a  subset $F$ of a poset $(P,\leq)$ is  a \textit{filter} if
$x\in F$ and $y\geq x$ imply $y\in F$. A
filter $F$ is {\it generated} by $X\subset P$ if $F=\langle
X\rangle=\{y\in P\,|\, y\geq x$ for some $x\in X\}$. If $x\in P$ then $\langle\{ x\}\rangle=\langle x\rangle$
is  a
\textit{principal filter}.  

The following two definitions quote   ~\cite[Definitions 1.5, 2.2,  and  3.1]{kt}. 

\begin{defn}\label{tight} The {\it hull} of a filter $F$ of $(\Gamma_0,\La)$ is the smallest filter $H_{\La}(F)$ of $(\Gamma_0,\La)$ containing $F$ and each vertex of $\Gamma_0\setminus F$ connected by an edge to a vertex in $F$.  A nonempty $S\in\mathfrak S$ is {\it principal}
 if  its
canonical form $S_{1}\dots S_r$ satisfies
$\mathrm{Supp}\,{S_{i}}=H_{\La}(\mathrm{Supp}\,{S_{i+1}})$ for
$0<i<r$ where $\mathrm{Supp}\,{S_{r}}$ is a principal filter.   The set $\mathfrak P$ of principal sequences inherits the binary relations $\sim$ and $\preccurlyeq$ from
$\mathfrak S$. By ~\cite[Proposition 1.9(b)]{kt}, a principal sequence is
determined uniquely up to equivalence by $r$ and
$\mathrm{Supp}\,{S_{r}}$, so we denote by $S_{r,x}$  the principal
sequence of size $r$ with $\mathrm{Supp}\,{S_{r}}=\langle x\rangle,\
x\in\Gamma_{0}$. 
\end{defn}

\begin{defn}If $S\in \mathfrak S$ annihilates a $k\gl$-module $M$, but no
proper subsequence of $S$ annihilates $M$, we call $S$ a {\it shortest} sequence annihilating $M$.
\end{defn}

The following statement quotes \cite[Proposition 2.1 and Theorems 2.2 and
2.6]{kt2}.

\begin{thm}\label{shrtstsq} Let $M$ be a preprojective $k\gl$-module.
\begin{itemize}
\item[(a)]  There exists a unique up to equivalence shortest  sequence
$S_M\in\mathfrak S$ annihilating $M$.
\item[(b)]  If $M$ is indecomposable and $N$ is an
indecomposable preprojective $k\gl$-module, then $S_N\sim S_M$ if and only if $N\cong M$.
\item[(c)]  If $M$ is indecomposable, then $S_M\in\mathfrak P$.
\item[(d)]  If $M$ is indecomposable and $S_M=x_1,\dots ,x_s$, then $M\cong F^-_{x_1}\dots F^-_{x_{s-1}}(L_{x_s})$ where $L_{x_s}$ is the simple projective $k(\Gamma,\sigma_{x_{s-1}}\dots\sigma_{x_{1}}\La)$-module associated with $x_s\in\Gamma_0$.
\end{itemize}
\end{thm}

\section{Reduced words in the Weyl group }\label{Section:Weylgroup}

For a graph $\Gamma$ with a valuation $\mathbf b$ we assume for the rest of the paper that $\Gamma_0=\{1,\dots,n\}$. Then the matrix $A=(a_{ij})$ with  $a_{ii} = 2$  and $a_{ij} = -b_{ij}$ for all $i\ne j$  is an indecomposable symmetrizable generalized $n\times n$
Cartan matrix, and $\mathcal W={\mathcal W}(A)$  is the Weyl group of the valued graph $(\Gamma,\mathbf b)$ ~\cite{dr}.    As before, $\gl$ denotes a quiver with a valuation $\mathbf b$ and modulation $\mathfrak B$.

\begin{defn}
If $S=x_{1},\dots,x_{s}$ is in $\mathfrak S$, we say that
$w(S)=\sigma_{x_{s}}\dots  \sigma_{x_{1}}$  is
the word in the Weyl group $\mathcal W$ associated to $S$.  If no
edge connects vertices $i$ and $j$, then
$\sigma_i\sigma_j=\sigma_j\sigma_i$ so that $S\sim T$ implies
$w(S)=w(T)$.
\end{defn}

Recall (see ~\cite{bourbLie4-6}) that the
{\it length} of $w\in\mathcal W$  is the smallest integer $\ell(w)=l\ge0$ such
that $w$ is the product of $l$ simple reflections, and a word
$w=\sigma_{y_{t}}\dots  \sigma_{y_{1}}$ in $\mathcal W$ is {\it
reduced} if $\ell(w)=t$.

We recall definitions and facts from \cite{ars,r}. Let ${\mathbb N}$ be the set of nonnegative integers.  The 
translation quiver ${\mathbb N}(\Gamma,\Lambda^{op})$ of the opposite quiver of 
$\gl$ has  ${\mathbb N}\times\Gamma_{0}$ as the set of vertices, and each 
arrow $a:u\to v$ of $\gl$  gives rise to two series of arrows, 
$(n,a^{\circ}):(n,v)\to(n,u)$ and 
$(n,a^{\circ})':(n,u)\to(n+1,v)$.  Since ${\mathbb 
N}(\Gamma,\Lambda^{op})$ is a locally finite quiver without oriented 
cycles (we disregard the valuation and translation), ${\mathbb N}\times\Gamma_{0}$ is a poset.

If
$X\in\mathrm{f.d.}\,k\gl$ is indecomposable, let $[X]$ be the
isomorphism class of $X$. If $Y\in\mathrm{f.d.}\,k\gl$ is
indecomposable, a path of length $m>0$ from $X$ to $Y$ is a sequence
of nonzero nonisomorphisms $X=A_{0}\to\dots\to A_{m}=Y,$ where
$A_{i}\in\mathrm{f.d.}\,k\gl$ is indecomposable for all $i$. By
definition, there is a path of length $0$ from $X$ to $X$.  Set $[X]\prec [Y]$ if there is a path of positive length from
$X$ to $Y$.

The preprojective component of $\gl$, $\tilde{\mathscr P}\gl$, is a
locally finite connected quiver  (we disregard the  valuation and translation) whose set of
vertices, $\tilde{\mathscr P}\gl_{0}$, is the set of isomorphism
classes of indecomposable preprojective $k\gl$-modules. If
$X,Y\in\mathrm{f.d.}\,k\gl$ are indecomposable, $Y$ is preprojective, and
$X=A_{0}\to\dots\to A_{m}=Y,\ m>0,$ is a path from $X$ to $Y$, then
$[X]\ne[Y]$ and $A_{i}$ is preprojective for all $i$.  Hence the reflexive closure $\preccurlyeq$ of the transitive binary
relation $\prec$ is a partial order on $\tilde{\mathscr P}\gl_{0}$.
Moreover, $[X]\prec [Y]$ if and only if there is a finite sequence
of irreducible morphisms $X=B_{0}\to\dots\to B_{n}=Y$, where $n>0$
and $B_{j}$ is indecomposable preprojective for all $j$.

The following two statements extend ~\cite[Proposition 3.7 and Theorem 4.3]{kp} to representations of valued quivers.
\begin{prop}\label{fullsubposet} If $M,N$ are
indecomposable preprojective $k\gl$-modules, then $[M]\preccurlyeq[N]$  in $\tilde{\mathscr P}\gl_{0}$ if and only if $S_M\preccurlyeq S_N$ in $\mathfrak P$.
\end{prop}
\begin{proof} The necessity follows from ~\cite[Corollary 2.9]{kt2}.  For the sufficiency, let  $S_M\preccurlyeq S_N$.   If $S_N\preccurlyeq S_M$, then $S_M\sim S_N$ by Proposition \ref{subsequence}(b) whence $M\cong N$ and $[M]=[N]$ according to Theorem \ref{shrtstsq}(b).  If $S_N\not\preccurlyeq S_M$ then, by Theorem \ref{shrtstsq}(a), $S_M\sim S_{p,u}$ and $S_N\sim S_{q,v}$ for some $p,q>0$ and $u,v\in\Gamma_0$.  By ~\cite[Theorem 2.5(a)]{kt}, which deals with $\mathfrak S$ but not with $\tilde{\mathscr P}$,   $(p-1,u)<(q-1,v)$ in ${\mathbb N}\times\Gamma_{0}$ so there is a path $(p-1,u)\to(q-1,v)$ of positive length in ${\mathbb N}(\Gamma,\Lambda^{op})$.  By ~\cite[Proposition 2.8(d)]{kt2}, there is a path $[M]\to[N]$ of positive length in $\tilde{\mathscr P}\gl$, i.e., $[M]\preccurlyeq[N]$. 
\end{proof} 

\begin{thm}\label{reducedword} Let $M$ be a preprojective $k\gl$-module.
\begin{itemize}
\item[(a)]  The word $w(S_M)\in \mathcal W$ is reduced.
\item[(b)]  If $M$ is
indecomposable and $N$ is an
indecomposable preprojective $k\gl$-module, the following are equivalent.
\begin{enumerate}
\item[(i)] $M\cong N$.
\item[(ii)] $S_M\sim S_N$.
\item[(iii)] $w(S_M)=w(S_N)$.
\end{enumerate}

\end{itemize}
\end{thm}
\begin{proof}
(a) If $M=0$ the statement is trivial. If $M\ne0$, let
$S_M=x_1,\dots,x_s$ and proceed by induction on $s>0$. The case
$s=1$ is clear, so suppose $s>1$ and the statement holds for all sequences $S_N$ of length $<s$.  Since $s>1$, $F^+_{x_1}M\ne0$ is
a preprojective $k(\Gamma,\sigma_{x_1}\La)$-module, and
$S_{F^+_{x_1}M}=x_2,\dots,x_s$.  By the induction hypothesis, the
word $u=\sigma_{x_s}\dots\sigma_{x_2}$ in $\mathcal W$ is reduced.
Assume, to the contrary, that the word
$u\sigma_{x_1}=\sigma_{x_s}\dots\sigma_{x_2}\sigma_{x_1}$ is not
reduced.  Then $\ell(u\sigma_{x_1})<\ell(u)$ \cite[Ch. IV,
Proposition 1.5.4]{bourbLie4-6} and
$\sigma_{x_s}\dots\sigma_{x_2}(e_{x_1})<0$ where $e_{x_1}$ is the
simple root associated to the vertex $x_1$ \cite[Lemma
3.11]{kac1990}.  By  ~\cite[Proposition 2.1]{dr},
$F(S_{F^+_{x_1}M})L_{x_1}=0$ where $L_{x_1}$ is the simple
$k(\Gamma,\sigma_{x_1}\La)$-module associated to $x_1$.  Since $x_1$ is a source in
$(\Gamma,\sigma_{x_1}\La)$, then $L_{x_1}$ is a simple injective and a
preprojective $k(\Gamma,\sigma_{x_1}\La)$-module so that
$[L_{x_1}]$ is a sink in $\tilde{\mathscr
P}(\Gamma,\sigma_{x_1}\La)$ and, hence, a
maximal element of the poset $\tilde{\mathscr
P}(\Gamma,\sigma_{x_1}\La)_0$. On the other hand, $F(S_{F^+_{x_1}M})L_{x_1}=0$  implies
$S_{L_{x_1}} \preccurlyeq S_{F^+_{x_1}M}$  whence $S_{L_{x_1}}\preccurlyeq S_{N}$ for some indecomposable direct
summand $N$ of $F^+_{x_1}M$ according to
\cite[Theorem 3.4(b)]{kp}, whose proof works for representations of
valued quivers.  By Proposition \ref{fullsubposet},
$[L_{x_1}]\preccurlyeq [N]$ in $\tilde{\mathscr
P}(\Gamma,\sigma_{x_1}\La)_0$.  Since
$[L_{x_1}]$ is a maximal element, $[L_{x_1}]=[N]$ whence
$L_{x_1}\cong N$.  The latter contradicts the well-known fact, contained in  ~\cite[Proposition 2.1]{dr}, that the simple module associated to a vertex that is
a source is not a direct summand of a module that belongs to the
image of the positive reflection functor associated to the vertex.
Thus $w(S_M)$ is a reduced word.

(b) The proof is identical to that of \cite[Theorem
4.3(b)]{kp}.
\end{proof}

With the same proof, the following statement extends \cite[Corollary 4.4]{kp} suggested by Zelevinsky.

\begin{cor}\label{zel}  Let $S=x_1,\dots,x_s,\ s>0,$ be in $\mathfrak S$, and set $M(S)=$ \newline$ F^-_{x_1}\dots  F^-_{x_{s-1}}(L_{x_s})$, where $L_{x_s}$ is the simple projective $k(\Gamma,\sigma_{s-1}\dots\sigma_{x_1}\La)$-module associated to $x_s\in\Gamma_0$.
\begin{itemize}
\item[(a)]  If the word $w(S)\in\mathcal W$ is reduced, $M(S)$ is an indecomposable module in $\tilde{\mathscr P}$.
\item[(b)]  If $M\in\tilde{\mathscr P}$ is indecomposable, then $M\cong M(S)$ for some sequence $S\in\mathfrak S$ where $\ell(S)>0$ and the word $w(S)$ is reduced.
\end{itemize}
\end{cor}

\begin{lem}\label{SeqSubShortest}
Let  $S = x_1, x_2, \dots, x_s,\ s >
1,$ be in $\mathfrak S$, suppose that the full subgraph of $\Gamma$ determined by $\mathrm{Supp}\,S$ is connected, and set $T = x_2 ,\dots, x_s$.
If $T \sim S_N$ for some indecomposable preprojective $k(\Gamma,\sigma_{x_1}\La)$-module $N$ satisfying
$M=F^-_{x_1}N \neq 0$, then $M$  is  indecomposable preprojective and $S \sim S_M$.
\end{lem}

\begin{proof} By  ~\cite[Proposition 2.1]{dr}, $M$ is  indecomposable and $F^+_{x_1}F^-_{x_1}N\cong N$.  Hence $S$ annihilates $M$ and $S_M
\preccurlyeq x_1 \vee S_M \sim x_1 U \preccurlyeq S = x_1 T$ for
some $U$, so by Proposition \ref{subsequence}(c), $U \preccurlyeq T$.  On the other hand, $x_1 U$
annihilates $M$, so $U$ annihilates $N$, giving $T \preccurlyeq U$,
hence $U \sim T$ and $x_1 \vee S_M \sim S$.  Since the full subgraph
of $\Gamma$ determined by $\mathrm{Supp}\, S$ is connected, for
some $x_l \in \mathrm{Supp}\, S \setminus \{ x_1 \}$ there is an
arrow $x_l \to x_1$ in $\gl$.  This gives   $x_l \in
\mathrm{Supp}\,S_M$ whence $x_1 \in
\mathrm{Supp}\,S_M$ because $\mathrm{Supp}\,S_M$ is  a filter of $(\Gamma_0,\La)$ by ~\cite[Proposition 1.9(a)]{kt}.  Hence $S_M \sim x_1 \vee S_M \sim S$.
\end{proof}

\begin{thm}\label{princseqreducedwrd}  If $S=x_1,\dots,x_s$, $s>0,$ is in $\mathfrak P$, the following are equivalent.
\begin{itemize}
\item[(a)]   There exists an indecomposable preprojective $k\gl$-module $M$ satisfying $S\sim S_M$.
\item[(b)]   The word $w(S)\in \mathcal W$ is reduced.
\item[(c)]  For $0<i<s,\ \sigma_{x_i}\dots\sigma_{x_{s-1}}(e_{x_s})>0$.
\end{itemize}
\end{thm}
\begin{proof} (a)$\implies$(b)\hskip.05in  This is Theorem \ref{reducedword}(a).

(b)$\implies$(c)\hskip.05in  This follows from Corollary
\ref{zel}(a).

(c)$\implies$(a)\hskip.05in   Proceed by  induction on $s$.  If $s=1$ then $S\sim S_{L_{x_1}}$ where $L_{x_1}$ is the simple projective $k\gl$-module associated to $x_1\in\Gamma_0$.  Suppose  $s>1$ and  the statement  holds for all principal (+)-admissible sequences of length $<s$ on all valued quivers $(\Gamma,\Theta)$ without oriented cycles.  By ~\cite[Proposition 3.6]{kp}, $T$ is a principal (+)-admissible sequence on $(\Gamma,\sigma_{x_1}\La)$, so by the induction hypothesis, $T\sim S_N$ for some indecomposable preprojective $k(\Gamma,\sigma_{x_1}\La)$-module $N$, and Theorem \ref{shrtstsq}(d) says that $N\cong F^-_{x_2}\dots
F^-_{x_{s-1}}(L_{x_s})$ where $L_{x_s}$ is the simple projective
$k(\Gamma,\sigma_{s-1}\dots\sigma_{x_1}\La)$-module
associated to $x_s\in\Gamma_0$.  Since $\sigma_{x_1}\dots\sigma_{x_{s-1}}(e_{x_s})>0$,  ~\cite[Proposition 2.1]{dr} implies that $M=F^-_{x_1}N\ne0$. By ~\cite[Remark 3.1]{kp}, the full subgraph of $\Gamma$ determined by
$\mathrm{Supp}\, S$ is connected.  By Lemma \ref{SeqSubShortest}, $M$ is indecomposable and $S
\sim S_M$.
\end{proof}

The proof of the following statement is identical to that of \cite[Theorem
4.6]{kp}.

\begin{thm}\label{seqreducedwrd} For all $S\in\mathfrak S$, the following are equivalent.
\begin{itemize}
\item[(a)]  There exists a preprojective $k\gl$-module $M$ satisfying $S\sim S_M$.
\item[(b)]  The word $w(S)\in \mathcal W$ is reduced.
\end{itemize}
\end{thm}

\begin{cor}\label{reduced} If    a valued graph $(\Gamma,\mathbf b)$  is not a Dynkin diagram of the type $A_n, B_n, C_n, D_n, E_6, E_7, E_8, F_4,$ or $G_2$, let $\La$ be an orientation on $\Gamma$ and  $S\in\mathfrak S$.
\begin{itemize}
\item[(a)] The word $w(S)  \in \mathcal W$ is reduced.
\item[(b)] For any  modulation $\mathfrak B$ of $(\Gamma,\mathbf b)$, there exists a preprojective $k\gl$-module $M$ satisfying $S\sim S_M$.
\end{itemize}
\end{cor}
\begin{proof}  (a) By ~\cite[Proposition, p. 224]{r2}, there exists a modulation $\mathfrak B$ of $(\Gamma,\mathbf b)$.  For this modulation the algebra $k\gl$ is of infinite representation
type by ~\cite[Theorem, p. 3]{dr}.  The poset $(\Gamma_0,\La)$ is finite,  so by Theorem
\ref{reducedword} there is an indecomposable preprojective
$k\gl$-module $M$ with $S_M \sim K^i U$ for some $i$ larger than the
size of $S$.  By Theorem \ref{seqreducedwrd}, $w(S_M) \in \mathcal W$
is reduced. Since $i$ is larger than the size of $S$, then $S_M \sim SV$
for some $V$, so $w(S) \in \mathcal W$ is also reduced.

(b) This is an immediate consequence of (a) and Theorem \ref{seqreducedwrd}.
\end{proof}

We characterize infinite Weyl groups in terms of reduced words.

\begin{thm}\label{infWeylGrp}
Let $A=(a_{ij})$ be an indecomposable symmetrizable generalized
$n\times n$ Cartan matrix, and let $c=\sigma_{v_n}\dots\sigma_{v_1}$
be a Coxeter element of the Weyl group $\mathcal W(A)$.   Then
$\mathcal W(A)$ is infinite if and only if for all $m\in\mathbb Z$,
$\ell(c^m)=|m|n$.
\end{thm}
\begin{proof}
Since sufficiency is clear, we need only prove necessity.  Consider the $n\times n$ matrix $(b_{ij})$ where $b_{ii}=0$ for all $i$, and $b_{ij}=-a_{ij}$ for all $i\ne j$.  Denote by $\Gamma=(\Gamma_0,\Gamma_1)$ the  graph where $\Gamma_0=\{1,\dots,n\}$ and,  for $i\ne j$, we have $\{i,j\}\in\Gamma_1$ if and only if $b_{ij}\ne0$. The collection of all $b_{ij}$ defines a  valuation $\mathbf b$ of $\Gamma$, and $\mathcal W(A)$ coincides with  the Weyl group $\mathcal W$ of the valued graph $(\Gamma,\mathbf b)$.  As noted
in \cite[p. 8]{dr}, there is a unique orientation $\La$ on $\Gamma$
 without oriented cycles for which 
$K=v_1,\dots,v_n$ is a complete sequence in $\mathfrak S$. Then
$c=w(K)$ and $c^m=w(K^m)$   for all $m>0$. Since $\mathcal W$ is
infinite, ~\cite[Propositions 1.2(a) and 1.5]{dr}  imply that $(\Gamma,\mathbf b)$ is not a
Dynkin diagram of the type $A_n, B_n, C_n, D_n, E_6, E_7, E_8, F_4,$ or
$G_2$. By Corollary \ref{reduced}(a), $c^m$ is a reduced word.
\end{proof}

\bibliographystyle{amsplain}

\begin{thebibliography}{10}
\bibitem[1]{ars}
M. Auslander, I. Reiten and S. O. Smal\o,
{\it Representation theory of Artin algebras},
Cambridge Studies in Advanced Mathematics, vol. 36,
Cambridge University Press, New York, 1994.

\bibitem[2]{bgp}
I. N. Bernstein, I. M. Gelfand, and V. A. Ponomarev, {\it Coxeter functors and
Gabriel's theorem}, Usp. Mat. Nauk {\bf 28} (1973), 19-33. Transl. Russ. Math.
Serv. {\bf 28} (1973), 17-32.

\bibitem[3]{bourbLie4-6}
N. Bourbaki, {\it Lie groups and Lie algebras. Chapters 4--6},
Translated from the 1968 French original by Andrew Pressley.
Elements of Mathematics (Berlin). Springer-Verlag, Berlin, 2002.
xii+300 pp.

\bibitem[4]{dr}
V. Dlab and C. M. Ringel, {\it Indecomposable representations of
graphs and algebras}, Memoirs Amer. Math. Soc., {\bf 173} (1976).

\bibitem[5]{dr1}
V. Dlab and C. M. Ringel, {\it On algebras of finite representation
type},  J. Algebra {\bf 33} (1975), 306--394.

\bibitem[6]{FZ2003}
S. Fomin and A. Zelevinsky, {\it $Y$-systems and generalized
associahedra.} Ann. of Math. (2) {\bf 158} (2003), no. 3, 977-1018.

\bibitem[7]{CAIV}
S. Fomin and A. Zelevinsky, {\it Cluster algebras IV: coefficients},
arXiv:math.RA/0602259.

\bibitem[8]{h}
R. B. Howlett, {\it Coxeter groups and $M$-matrices}. Bull. London
Math. Soc. {\bf 14} (1982), no. 2, 137--141.

\bibitem[9]{kac1990}
V. G. Kac, {\it Infinite-dimensional Lie algebras. Third edition},
Cambridge University Press, Cambridge, 1990. xxii+400 pp.

\bibitem[10]{kp}
M. Kleiner and A. Pelley, {\it Admissible sequences, preprojective
modules, and reduced words in the Weyl group of a quiver},
arXiv:math.RT/0607001.

\bibitem[11]{kt}
M. Kleiner and H. R. Tyler, {\it Admissible sequences and the
preprojective component of a quiver}, Adv. Math. {\bf 192} (2005),
no. 2, 376--402.

\bibitem[12]{kt2}
M. Kleiner and H. R. Tyler, {\it Sequences of reflection functors and the preprojective component of a valued quiver}, arXiv:math.RT/0608175.

\bibitem[13]{mrz2003}
R. Marsh, M.  Reineke, and A. Zelevinsky, {\it Generalized
associahedra via quiver representations.} Trans. Amer. Math. Soc.
{\bf 355} (2003), no. 10, 4171--4186.

\bibitem[14]{McCam}
J. McCammond, {\it Noncrossing partitions in surprising locations.}
Amer. Math. Monthly, to appear.

\bibitem[15]{rea}
N. Reading, {\it Clusters, Coxeter-sortable elements and noncrossing
partitions}, Trans. Amer. Math. Soc., in press,
arXiv:math.CO/0507186.

\bibitem[16]{r2}
C. M. Ringel,  \emph{Green's theorem on {H}all algebras},
Representation theory of algebras and related topics (Mexico City, 1994), CMS Conf.
Proc.,
  vol.~19, Amer. Math. Soc., Providence, RI, 1996, pp.~185--245.

\bibitem[17]{r}
C. M. Ringel, {\it Tame algebras and integral quadratic forms},
Lecture Notes in Math., vol. 1099, Springer-Verlag, Berlin, 1984.
\end{thebibliography}

\end{document}